\documentclass{amsart}

\usepackage{amssymb}

\newtheorem{theorem}{Theorem}[section]

\newtheorem{corollary}[theorem]{Corollary}
\def\square#1{\vbox{\hrule
\hbox{\vrule\hbox to #1 pt{\hfill}\vbox{\vskip #1 pt}\vrule}\hrule}}

\theoremstyle{definition}
\newtheorem{definition}[theorem]{Definition}
\newtheorem{example}[theorem]{Example}

\theoremstyle{remark}
\newtheorem{remark}[theorem]{Remark}

\numberwithin{equation}{section}

\font\a=msbm10
\font\aa=msbm10 at 8pt
\font\d=cmr10

\font\eee=cmcsc10 at 8pt
\font\bbb=cmbx10 at 8pt

\begin{document}

\title{Multiplicity of complex hypersurface singularities, Rouch\'e satellites and
  Zariski's problem} 

\author{Christophe Eyral}
\address{Max-Planck Institut f\"ur Mathematik, Vivatsgasse 7, 53111
Bonn, Germany}  
\email{eyralchr@yahoo.com}

\author{Elizabeth Gasparim}
\address{Department of Mathematical Sciences, New Mexico State
University, Las Cruces NM 88003-8001, USA} 
\email{gasparim@nmsu.edu}



\thanks{This research was supported by the Max--Planck Institut
 f\"ur Mathematik in Bonn.} 

\begin{abstract}
Let $f,g\colon\, (\hbox{\aa C}^n,0) \rightarrow (\hbox{\aa C},0)$ be 
reduced germs  of holomorphic functions.  
We show that $f$ and $g$ have the same multiplicity at $0$, if and
only if, there exist reduced germs $f'$ and $g'$ analytically
equivalent to  $f$ and $g$, respectively, such that $f'$ and $g'$  satisfy a
Rouch\'e type inequality with respect to a generic  
`small' circle around $0$. As an application, we give a
reformulation of Zariski's multiplicity question and a partial positive 
answer to it.

\medskip
\noindent{\eee R\'esum\'e.} {\bbb Multiplicit\'e des singularit\'es
  d'hypersurfaces complexes, satellites de Rouch\'e et probl\`eme de
  Zariski.} Soient $f,g\colon\, (\hbox{\aa C}^n,0) 
\rightarrow (\hbox{\aa C},0)$ des germes de fonctions
holomorphes r\'eduits. Nous montrons que $f$ et $g$ ont 
la m\^eme multiplicit\'e
en $0$ si et seulement s'il existe des germes r\'eduits $f'$ et
$g'$ analytiquement \'equivalents \`a $f$ et $g$, respectivement, tels
que $f'$ et $g'$ satisfassent une in\'egalit\'e du type de Rouch\'e
par rapport \`a un `petit' cercle g\'en\'erique autour de~$0$. Comme
application, nous donnons une reformulation de la question de Zariski
sur la multiplicit\'e et une r\'eponse partielle positive \`a celle--ci.
\end{abstract}

\maketitle

\markboth{Christophe Eyral and Elizabeth Gasparim}{Multiplicity of
  complex hypersurface singularities, 
  Rouch\'e satellites and 
  Zariski's problem}

\section{Introduction}

Let $f,g\colon\, (\hbox{\a C}^n,0) \rightarrow (\hbox{\a C},0)$ be reduced
germs (at the origin) of  holomorphic functions, with $n\geq 2$,
$V_f$, $V_g$ the corresponding germs of
hypersurfaces in $\hbox{\a C}^n$, and $\nu_f, \nu_g$ the
multiplicities at~$0$ of $V_f$, $V_g$ respectively.
By the {\it multiplicity}\/ $\nu_f$ we mean
the number of points of intersection, near~$0$,
of~$V_{f}$ with a generic (complex) line in 
$\hbox{\a C}^n$ passing arbitrarily close to $0$ but not
through~$0$. As we are assuming that $f$ is reduced, 
$\nu_{f}$ is also the {\it order} of $f$ at $0$, that is, the lowest degree
in the power series expansion of~$f$ at~$0$.  We denote by $C(V_f)$,
$C(V_g)$ the tangent cones at $0$ of $V_f$, 
$V_g$, that is, the zero sets of the initial polynomials
of $f$ and $g$ respectively (cf.~\cite{W}).

In Section 2, we prove that $\nu_f=\nu_g$, if and
only if, there exist reduced germs $f'$ and $g'$ analytically equivalent to $f$ and
$g$, respectively, such that $\vert f'(z)-g'(z)\vert < \vert f'(z)\vert$,
for all $z\in \dot D$, where $\dot D$ is the boundary of a generic
`small' disc around $0$ (Theorem~\ref{t:t3}). We call such an
inequality a {\it Rouch\'e inequality}\/ and we say that $g'$ is a
{\it Rouch\'e satellite}\/ of $f'$.

In Section 3, we apply this result to Zariski's
multiplicity question. In particular, we show that the answer to 
Zariski's question is {\it yes}, if and only if, for any two
topologically equivalent reduced germs $f$ and $g$ there exist reduced
germs $f'$ and $g'$
analytically equivalent to $f$ and $g$, respectively, such that
$g'$ is a Rouch\'e satellite of $f'$ (Theorem \ref{t:tz}). In addition, we
answer positively Zariski's question in the special 
case of `small' homeomorphisms for Newton nondegenerate isolated 
singularities (Corollary \ref{t:c2}) and one--parameter families of isolated 
singularities (Corollary \ref{families}).

\section{Multiplicity and Rouch\'e satellites}

Let $L$ be a line through $0$ in $\hbox{\a C}^n$ not contained
in $C(V_f)\cup C(V_g)$ (equivalently, $L\cap
(C(V_f)\cup C(V_g))=\{0\}$). Then~$\nu_f$ (respectively $\nu_g$) is the
order at $0$ of $f_{\mid L}$ (respectively~$g_{\mid L}$), and $0$ is
an isolated point of $L\cap V_f$ and $L\cap V_g$ (cf.~\cite{Ep}). 
In particular, there exists a closed 
disc~$D\subseteq L$ around~$0$ such that, for any closed disc
$D'\subseteq D$ around~$0$, $D'\cap (V_f\cup V_g)=\{0\}$.
We shall call such a disc $D$ a {\it good}\/ disc
for $f$ and for $g$. 

\begin{definition}\label{satellite}
We say that $g$ is a {\it Rouch\'e satellite}\/ of $f$ if there exists a
good disc $D$ (for $f$ and for $g$) such that $f$ and $g$ satisfy a
{\it Rouch\'e inequality}\/ with respect to the boundary $\dot D$ of $D$,
that is, $$\vert f(z)-g(z)\vert < \vert f(z)\vert$$
for all $z\in \dot D$. 
\end{definition}

\begin{theorem}\label{t:t1} 
If $g$ is a Rouch\'e satellite of $f$, then $\nu_g=\nu_f$.
\end{theorem}

\begin{proof}
Let $D\subseteq L$ be a good disc for $f$ and for $g$ (for some
line $L$ through $0$ not contained in $C(V_f)\cup C(V_g)$)
such that $\vert f_{\mid L}(z)-g_{\mid L}(z)\vert < \vert f_{\mid
  L}(z)\vert$ for all $z\in \dot D$.
By Rouch\'e theorem (cf.~e.g.~\cite[Chapter VI, Theorem 1.6]{L}),
$f_{\mid L}$ and $g_{\mid L}$ have the same number of zeros, counted with their
multiplicities, in the interior of $D$. Thus, since $f_{\mid L}$ and $g_{\mid L}$
vanish only at $0$ on~$D$, the orders at $0$ of $f_{\mid L}$ and $g_{\mid L}$
are equal. In other words, $\nu_f=\nu_g$.
\end{proof}

\begin{example}\label{ex1}
Consider the germs $f,\, g\colon\, (\hbox{\a C}^3,0) \rightarrow
(\hbox{\a C},0)$ defined by
$$f(z_1,z_2,z_3)=z_1^2+z_2^3+z_3^3+z_1^3+z_2^4\quad \hbox{and} \quad
g(z_1,z_2,z_3)=z_1^2+z_2^3+z_3^3+z_1^4+z_2^6.$$
Then $g$ is a Rouch\'e satellite of $f$. Indeed,  
set $L=\{(z_1,0,z_3)\in\hbox{\a C}^3 \mid z_1=z_3\};$  
then 
 $$V_f\cap L=\Bigl\{(0,0,0),\Bigl(-\frac{1}{2},0,-\frac{1}
{2}\Bigr)\Bigr\}\quad\hbox{and}\quad  V_g\cap
L=\{(0,0,0),(a,0,a),(\bar a,0, \bar a)\},$$ where $a=(-1-i\sqrt{3})/2$
and $\bar a$ is the complex conjugate of $a$.
So, the disc  $D \subseteq L$  of radius
$1/4$ is good for $f $ and for $g,$ 
and, for all $z\in \dot D$, $$\vert f(z)-g(z)\vert\leq \frac{5}{4^4} <
\frac{2}{4^3}\leq \vert f(z)\vert.$$ 
Hence $g$ is a Rouch\'e satellite of $f$.
In fact, here, $f$ is also a Rouch\'e satellite
of $g$. Indeed, for all $z\in\dot D$, we have
$$\vert f(z)-g(z)\vert\leq \frac{5}{4^4} <
\frac{11}{4^4}\leq \vert g(z)\vert.$$ 
\end{example}

Of course, in general, $g$ may be a Rouch\'e satellite of
$f$ without  $f$ being a Rouch\'e satellite of $g$. For example,
take $g=f/2$. Also, it is not difficult to construct $f$ and $g$ such that
$\nu_f=\nu_g$ but neither $g$ is a Rouch\'e satellite of $f$ nor $f$  a
Rouch\'e satellite of~$g$. Take for example $g=-f$. Nevertheless,
such an unpleasant situation is resolved by Theorem \ref{converse} below.

\begin{definition}
If there exists a germ of homeomorphism 
$\varphi\colon (\hbox{\a C}^n,0)\rightarrow (\hbox{\a C}^n,0)$ such that:
\begin{enumerate} 
\item $\varphi(V_g)=V_f$ then $f$ and $g$ 
are called  {\it topologically  equivalent} (denoted $f\sim_t g);$
\item $\varphi(V_g)=V_f$ and $\varphi$ is an analytic isomorphism, 
then $f$ and $g$  
are called  {\it analytically  equivalent} (denoted $f\sim_a g);$
\item  $g= f\circ \varphi$ then $f$ and $g$ 
are called  {\it topologically right equivalent} (denoted $f\sim_{tr} g).$
\end{enumerate}
\end{definition}

Note that the definition makes sense only for {\it reduced}\/
germs. In the special case of an isolated singularity, the hypothesis
`$n\geq 2$' automatically implies that the germ is reduced. Note also
that $(2) \Rightarrow (1)$ and $(3) \Rightarrow (1).$

Theorem \ref{t:t1} has the weak following converse.

\begin{theorem}\label{converse}
If $\nu_f=\nu_g$, then there exist reduced germs $f'\sim_a f$ and $g'\sim_a g$
such that $g'$ is a Rouch\'e satellite of $f'$. 
\end{theorem}

\begin{proof}
By an analytic change of coordinates, one can assume that the
$z_n$--axis, $Oz_n$, is not con\-tained in the tangent cones $C(V_f)$,
$C(V_g)$, so that  
$f(0,\ldots,0,z_n)\not=0$ and $g(0,\ldots,0,z_n)\not=0$, for any $z_n\not=0$
close enough to $0$. By the Weierstrass preparation theorem, for $z$ 
near~$0$, the germ $f(z)$  can be represented as a
product $f(z)=f'(z)\, f''(z)$, where $f''(z)$ is a germ of holomorphic function
which does not vanish around $0$ and where $f'(z)$ is of the form
$$f'(z_1,\ldots,z_n)=z_n^{\nu_f} + z_n^{\nu_f-1}
f_1(z_1,\ldots,z_{n-1}) +\ldots + f_{\nu_f}(z_1,\ldots,z_{n-1}),$$
with, for $1\leq i \leq \nu_f$, $f_i\in\hbox{\a
  C}\{z_1,\ldots,z_{n-1}\}$, $f_i(0)=0$ and the order of $f_i$ at $0$
is $\geq i$. Similarly $g(z)=g'(z)\, g''(z)$, with $g''(z)\not=0$ for
all $z$ near $0$, and   
$$g'(z_1,\ldots,z_n)=z_n^{\nu_g} + z_n^{\nu_g-1}
g_1(z_1,\ldots,z_{n-1}) +\ldots + g_{\nu_g}(z_1,\ldots,z_{n-1}),$$
with, for $1\leq i \leq \nu_g$, $g_i\in\hbox{\a
  C}\{z_1,\ldots,z_{n-1}\}$, $g_i(0)=0$ and the order of $g_i$ at $0$
is $\geq i$. Clearly $f'$ and $g'$ are reduced, and,
since $V_f=V_{f'}$ and  $V_g=V_{g'}$, $f'\sim_a f$ and $g'\sim_a
g$. On the other hand, since $\nu_f=\nu_g$,
$f'_{\mid Oz_n}=g'_{\mid Oz_n}$. 
But for any disc $D\subseteq Oz_n$ around $0$ (in particular for
any good disc in $Oz_n$
for $f'$ and $g'$), $\vert f'(z)\vert=r^{\nu_f}\not=0$ for all $z\in\dot 
D$, where $r$ is the radius of~$D$.
\end{proof}

Since the multiplicity is an invariant of the (embedded) reduced
analytic type, 
we can summarize Theorems \ref{t:t1} and \ref{converse} as follows. 

\begin{theorem}\label{t:t3}
The multiplicities $\nu_f$ and $\nu_g$ are the
same, if and only if, there exist reduced germs $f'\sim_a f$ and
$g'\sim_a g$ such 
that $g'$ is a Rouch\'e satellite of $f'$.
\end{theorem}

\section{Applications to Zariski's multiplicity question}\label{section4}

In \cite{Z}, Zariski posed the following question: {\it if $f \sim_t g,$ 
then is it true that $\nu_f=\nu_g$?}  
The question is, in general, still unsettled (even for
hypersurfaces with isolated singularities). The answer is,
nevertheless, known to be {\it yes}\/ in several special cases the
list of which can be found in the recent first author's 
survey article \cite{Ey}. In particular, Ephraim \cite{Ep} proved that
multiplicity is preserved by ambient $C^1$--diffeomorphisms; his paper
inspired some of our proofs. In this section, we give a partial positive
answer to Zariski's question in the special 
case of `small' homeomorphisms for Newton nondegenerate isolated 
singularities and one--parameter families of isolated 
singularities. In addition, we give an equivalent reformulation of
Zariski's question in terms of Rouch\'e satellites.

We start with the following result which asserts that if $f$ and $g$ are
topologically right equivalent via a sufficiently `small' homeomorphism,
then they have the same multiplicity. More precisely suppose
$f\sim_{tr}g$. Then there are representatives 
$\hbox{f}\colon U\rightarrow\hbox{\a C}$ and $\hbox{g}\colon
U'\subseteq U\rightarrow\hbox{\a C}$ of the germs $f$ and
$g$ respectively  and a homeomorphism
$\varphi \colon U'\rightarrow \varphi(U')\subseteq U$ such that
$\varphi(0)=0$ and $\hbox{g}=\hbox{f}\circ \varphi$. Since $\hbox{f}$ is
uniformly continuous  
on a compact small ball $B_r\subseteq U'$ around $0$, there
exists $\eta>0$ such that, for any $z,w\in B_r$,  
$$\vert z-w\vert < \eta \ \Rightarrow \
\vert \hbox{f}(z)-\hbox{f}(w)\vert < \inf_{u\in \dot D_\varrho} \vert
\hbox{f}(u)\vert,$$ 
where $D_\varrho$ is a good disc at $0$ for $\hbox{f}$ and for
$\hbox{g}=\hbox{f}\circ\varphi$ with radius $\varrho\leq r/2$.

\begin{definition}
We will say that the homeomorphism $\varphi \colon U'\rightarrow
\varphi(U')\subseteq U$ is $\hbox{f}${\it --small}\/ if there exists a triple
$(r,\varrho,\eta)$ as above 
such that, for all $z\in B_r$, $$\vert z-\varphi(z)\vert < \inf\{\eta,
\varrho\}.$$ 
\end{definition}

\begin{theorem}
With the above hypotheses and notation, if the homeomorphism $\varphi
\colon U'\rightarrow 
\varphi(U')\subseteq U$ is $\hbox{\d f}$--small, then $\nu_{f}=\nu_{g}$.   
\end{theorem}

\begin{proof}
By hypothesis, for all $z\in \dot D_{\varrho}$, $\varphi(z)\in
B_r$ and $$\vert \hbox{f}(z)-\hbox{f}\circ \varphi(z)\vert <
\inf_{u\in \dot D_\varrho} 
\vert \hbox{f}(u)\vert \leq \vert \hbox{f}(z)\vert.$$  
Therefore $\hbox{g}=\hbox{f}\circ \varphi$ is a Rouch\'e satellite of
$\hbox{f}$. Then, by Theorem \ref{t:t1}, $\nu_{f}=\nu_{g}$. 
\end{proof}

The interest in topologically right equivalent germs with regard to
Zariski's question comes from the following.
By theorems of King \cite{K}, Perron \cite{P}, Saeki
\cite{S} and Nishimura~\cite{N}, if $f$ has an 
{\it isolated}\/ singularity at $0$ and a nondegenerate Newton
principal part, then the relation $f \sim_t g$ implies 
$f \sim_{tr} g$. 
On the other hand, by another theorem of King \cite{Ki2}, for a
one--parameter holomorphic family  of {\it isolated}\/ 
singularities $(f_s)_s$ in $\hbox{\a C}^n$, with $n\not=3$, if the
relation $f_s 
\sim_t f_0$  holds for all $s$ near $0$,  
then so does $f_s \sim_{tr} f_0$. So, when considering
 isolated Newton nondegenerate 
singularities or {\it families}\/ of isolated singularities, the
Zariski problem refers immediately to right equivalent germs.

\begin{corollary}\label{t:c2}
Assume that $f$ has an isolated critical point at $0$ and a
nondegenerate Newton principal part, and suppose $g\sim_{t} f$. In
this case, there are representatives
$\hbox{\d f}\colon U\rightarrow\hbox{\a C}$ and $\hbox{\d g}\colon
U'\subseteq U\rightarrow\hbox{\a C}$ of $f$ and $g$ respectively and
a homeomorphism 
$\varphi \colon U'\rightarrow \varphi(U')\subseteq U$ such that
$\varphi(0)=0$ and $\hbox{\d g}=\hbox{\d f}\circ \varphi$. If 
$\varphi$ is $\hbox{\d f}$--small, then $\nu_f=\nu_g$. 
\end{corollary}

\begin{remark} If, in addition, $f$ is {\it convenient}\/ (cf.~
  \cite{Ko}), 
then the hypothesis of having an  
isolated singularity at $0$ is automatically satisfied (cf.~\cite{O2}).
\end{remark}

Corollary \ref{t:c2} is complementary to the result of Abderrahmane
and Saia--Tomazella  concerning
$\mu$--constant {\it families}\/ of convenient Newton nondegenerate (isolated)
singularities (cf.~\cite{Ab} and~\cite{ST}). 

\begin{corollary}\label{families}
Let $(f_s)_s$ be a topologically
constant (or $\mu$--constant) 
one--parameter holomorphic 
family of isolated hypersurface singularities, with $n\not=3$.
In this case, for all $s$ near~$0$, there are representatives
$\hbox{\d f}_0\colon U_0\rightarrow\hbox{\a C}$ and
$\hbox{\d f}_s\colon U_s\subseteq U_0\rightarrow\hbox{\a C}$ of $f_0$
and $f_s$ respectively and a
homeomorphism 
$\varphi_s \colon U_s\rightarrow \varphi(U_s)\subseteq U_0$ such that
$\varphi_s(0)=0$ and $\hbox{\d f}_s=\hbox{\d f}_0\circ \varphi_s$. If,
for all $s$ near~$0$, $\varphi_s$ is $\hbox{\d f}_0$--small, then $(f_s)_s$ is
equimultiple (i.e., for all $s$ near~$0$, $\nu_{f_s}=\nu_{f_0}$).  
\end{corollary}

We conclude with the following nice consequence of Theorem \ref{t:t3}
which is
reformulation of Zariski's multiplicity question in terms of
Rouch\'e satellites.

\begin{theorem}\label{t:tz}
The answer to Zariski's multiplicity question is yes,
if and only if,  the relation $f \sim_t g$ implies  that 
there  exist reduced germs $f' \sim_{a} f$ and $g'\sim_{a} g$ such
that $g'$ is a Rouch\'e satellite of $f'$.
\end{theorem}

\bibliographystyle{amsplain}

\begin{thebibliography}{10}

\bibitem {Ab} O.M. Abderrahmane,  \textit{On the deformation with
constant Milnor number and Newton polyhedron}, Preprint (Saitama
University, 2004).

\bibitem {Ep} R. Ephraim,  \textit{$C^1$ preservation of the multiplicity},
Duke Math.~J.~\textbf{43} (1976) 797--803.

\bibitem {Ey} C. Eyral,  \textit{Zariski's multiplicity question -- A survey},
Preprint (Max-Planck Institut f\"ur Mathematik, 2005).

\bibitem {K} H. King,   \textit{Topological type of isolated critical points},
Ann.~of Math.~\textbf{107} (1978) 385--397.

\bibitem {Ki2} H. King,   \textit{Topological type in families of germs},
Invent.~Math.~\textbf{62} (1980) 1--13.

\bibitem {Ko} A.G. Kouchnirenko,  \textit{Poly\`edres de Newton et nombres de
Milnor}, Invent.~Math.~\textbf{32} (1976) 1--32.

\bibitem {L} S. Lang,  \textit{Complex analysis}, Graduate Texts in
  Mathematics \textbf{103} (Springer--Verlag, corrected second
  printing, 1995).

\bibitem {P}  B. Perron,  \textit{Conjugaison topologique des germes de
fonctions holomorphes \`a singularit\'e isol\'ee en dimension trois}, 
Invent.~Math.~\textbf{82} (1985) 27--35.

\bibitem {N}  T. Nishimura,  \textit{A remark on topological types of complex
isolated singularities of  hypersurfaces}, private communication
between O. Saeki and T. Nishimura as cited in \cite {S}.

\bibitem {O2} M. Oka,  \textit{Non-degenerate complete intersection
  singularity}, Actualit\'es Math\'ematiques (Hermann, Paris, 1997).


\bibitem {S}  O. Saeki,  \textit{Topological types of complex isolated
hypersurface singularities}, Kodai Math.~J. \textbf{12} (1989)
23--29.

\bibitem {ST} M.J. Saia and J.N. Tomazella,  \textit{Deformations with
constant Milnor number and multiplicity of complex hypersurfaces},
Glasg.~Math.~J.~\textbf{46} (2004) 121--130.

\bibitem {W} H. Whitney, \textit{Complex analytic varieties}, Addison-Wesley
publishing company, Reading, Mass.-Menlo Park, Calif.-London-Don Mills, Ont,
1972.

\bibitem {Z} O. Zariski, \textit{Open questions in the theory of
singularities}, Bull.~Amer. Math.~Soc.~\textbf{77} (1971) 481--491.

\end{thebibliography}

\end{document}